\documentclass[%
 aip,
cp,  
 amsmath,amssymb,
 reprint,%
]{revtex4-2}

\usepackage{graphicx}
\usepackage{dcolumn}
\usepackage{bm}

\usepackage[utf8]{inputenc}
\usepackage[T1]{fontenc}
\usepackage{mathptmx} 

\begin{document}

\title{Some Aspects of the Numerical Analysis of a Fractional Duffing
	Oscillator with a Fractional Variable Order Derivative of the Riemann-Liouville Type}

\author{Kim Valentine}%
\affiliation{Vitus Bering Kamchatka State University, Petropavlovsk-Kamchatskiy, Russia 
}

\author{Parovik Roman} 
\email[Corresponding author: ]{romanparovik@gmail.com}
\affiliation{Vitus Bering Kamchatka State University, Petropavlovsk-Kamchatskiy, Russia 
}
\affiliation{Institute for Cosmophysical Research and Radio Wave Propagation, Far East Branch, Russian Academy of Sciences}


\begin{abstract}
In this paper, we consider some aspects of the numerical
analysis of the mathematical model of fractional Duffing with a derivative of
variable fractional order of the Riemann-Liouville type. Using numerical methods: an explicit finite-difference scheme based on the Grunwald-Letnikov and Adams-Bashford-Moulton approximations (predictor-corrector), the proposed
numerical model is found. These methods have been verified with a test case. It is shown that the predictor-corrector method has a faster convergence than the method according to the explicit finite-difference scheme. For these schemes, using Runge's rule, estimates of the computational accuracy were made, which tended to unity with an increase in the number of calculated grid nodes.
\end{abstract}

\maketitle

\section{Introduction}

Currently, one of the scientific directions of nonlinear dynamics, fractional
dynamics, has received wide development \cite{1}. She studies the hereditarity
properties of dynamical systems. Hereditarity (memory) is a property of a dynamic system in which its current state depends on its previous states \cite{2}. As shown in \cite{1}, the property of memory can be described using the mathematical apparatus of fractional calculus or using fractional derivative operators. Operators of fractional derivatives have many definitions and unique properties, but all of them, to one degree or another, describe the memory effect that characterizes information about the previous states of the system. This effect predetermines additional degrees of freedom -- orders of fractional derivatives \cite{3}. Such multi-parameter dynamical systems have certain difficulties in description and require special research methods to detect chaotic regimes \cite{4}. 

In this work, using the Adams-Bashford-Moulton method, a numerical solution of the Duffing equation with a derivative of variable fractional order was found. The fractional derivative operator was taken in the Riemann-Liouville sense. A comparison was made with the results of \cite{5}, in which the numerical solution was presented in the form of an explicit finite-difference scheme.

\section{MATHEMATICAL MODEL}

Consider the following Cauchy problem \cite{5}:

\begin{equation}
\ddot x\left( t \right) + \lambda D_{0t}^{q\left( t \right)}x\left( \tau 
\right) + \omega _0^2x\left( t \right) + b{x^3}\left( t \right) = f(t),\;x\left(
0 \right) = {x_0},\dot x\left( 0 \right) = {y_0},
\label{eq1}
\end{equation}
where $x\left( t \right) \in {C^2}\left( {0,T} \right)$ is the displacement
function, $\ddot x\left( t \right) = {{{d^2}x\left( t \right)} \mathord{\left/ {\vphantom {{{d^2}x\left( t \right)} {d{t^2}}}} \right.
\kern-\nulldelimiterspace} {d{t^2}}},\dot x\left( t \right) = dx/dt$, $\lambda
$ is the friction coefficient, $f(t) = \delta \cos (\omega t)$ is the external
action, $\delta $  and $\omega $ are the amplitude and frequency of the external periodic action, $\omega _0^2$  is the natural frequency of the system, $b$ is the isochronism coefficient, ${x_0}$ and ${y_0}$  are the given constants that determine the initial conditions, $0 < q\left( t \right) < 1$  is continuous function, $t$  -- time, $T$ -- simulation time.

The operator $D_{0t}^{q(t)}x\left( \tau  \right)$ in the model equation (1) has the form \cite{5}:
\begin{equation}
D_{0t}^{q(t)}x\left( \tau  \right) = \frac{1}{{\Gamma \left( {1 - q\left( t
			\right)} \right)}}\frac{d}{{dt}}\int\limits_0^t {\frac{{x\left( \tau 
			\right)d\tau }}{{{{\left( {t - \tau } \right)}^{q\left( t \right)}}}}}
\label{eq2}
\end{equation}

Operator (2) will be called the operator of the fractional variable order
derivative of the Riemann-Liouville type. Equation (1) is a fractional Duffing
equation. It is convenient to represent the Cauchy problem (1) in the form of a system of differential equations \cite{5}:
\begin{equation}
\left\{ \begin{array}{l}
\dot x\left( t \right) = y\left( t \right),\\
D_{0t}^{q(t)}x\left( \tau  \right) = w\left( t \right),\\
\dot y\left( t \right) = f(t) - \lambda w\left( t \right) - \omega _0^2x\left( t
\right) - b{x^3}\left( t \right),\\
x\left( 0 \right) = {x_0},\;\dot x\left( 0 \right) = {y_0}.
\end{array} \right.
\label{eq3}
\end{equation}

The variable order $q\left( t \right)$ of the fractional derivative of the
Riemann-Liouville type determines the intensity of energy dissipation in the
oscillatory system and is associated with the properties of the medium in which the oscillatory process takes place. In the case when this order is constant and equal to one, then the Cauchy problem (1) turns into the Cauchy problem for the classical Duffing oscillator.

\section{RESEARCH METHODOLOGY}

System (3), due to nonlinearity, does not have an exact solution; therefore, we will seek an approximate solution using the theory of finite-difference schemes. Divide the segment $\left[ {0,T} \right]$ into $N$  equal parts with step $h$. The solution to the differential problem $x\left( t \right)$ turns into an approximate grid solution $x\left( {{t_k}} \right),\;{t_k} = kh,\;k = 1,...,N$. The fractional derivative in system (3) is approximated by a difference analogue -- the Grunwald-Letnikov derivative \cite{5}:
\begin{equation}
D_{0t}^{q\left( t \right)}x\left( \tau  \right) \approx {w_{k - 1}} =
\frac{1}{{{h^{{q_k}}}}}\sum\limits_{j = 0}^{k - 1} {c_j^{\left( {{q_k}}
		\right)}{x_{k - j}}}  = \frac{{{x_k}}}{{{h^{{q_k}}}}} +
\frac{1}{{{h^{{q_k}}}}}\sum\limits_{j = 1}^{k - 1} {c_j^{\left( {{q_k}}
		\right)}{x_{k - j}}} ,
\label{eq4}
\end{equation}
\[
c_0^{\left( {{q_k}} \right)} = 1,c_j^{\left( {{q_k}} \right)} = \left( {1 -
	\frac{{1 + {q_k}}}{j}} \right)c_{j - 1}^{{q_k}}.
\]

Integer derivatives:
\begin{equation}
\dot x\left( t \right) \approx \frac{{{x_k} - {x_{k - 1}}}}{h},\dot y\left( t
\right) \approx \frac{{{y_k} - {y_{k - 1}}}}{h}
\label{eq5}
\end{equation}

Substituting (4) and (5) into system (3), we arrive at the following approximate solution of the Cauchy problem (1) \cite{5}:
\begin{equation}
\left\{ \begin{array}{l}
{x_k} = h{y_{k - 1}} + {x_{k - 1}}\\
{y_k} = {y_{k - 1}} + h\left( { - {x_{k - 1}} - bx_{k - 1}^3 + f({t_{k - 1}}) -
	\lambda {h^{ - {q_k}}}\sum\limits_{i = 0}^{k - 1} {c_i^{({q_k})}{x_{k - i}}} }
\right)\\
c_k^{({q_k})} = \left( {1 - \frac{{1 + {q_k}}}{k}} \right)c_{k -
	1}^{({q_k})},{c_0} = 1.
\end{array} \right.
\label{eq6}
\end{equation}

The discrete system (6) approximating the Cauchy problem (1) will be called the explicit finite difference scheme (EFDS). The results obtained by scheme (6) were considered in detail in \cite{5}.

Consider another method for solving the Cauchy problem (1) -- the
Adams-Bashfort-Moulton (ABM) method from the family of predictor-corrector
methods \cite{6}, \cite{7}. For this, we represent the Cauchy problem (1) in the form of
the system:
\begin{equation}
\left\{ \begin{array}{l}
\partial _{0t}^{{q_1}}x = y,{q_1} = 1\\
\partial _{0t}^{{q_2}(t)}y = z,{q_2}\left( t \right) = 1 - q(t)\\
\partial _{0t}^{{q_3}(t)}z = f(t) - \lambda z - x - b{x^3},{q_3}(t) = q(t).
\end{array} \right.
\label{eq7}
\end{equation}

We seek a solution to system (7) in the form:
\begin{equation}
\left\{ \begin{array}{l}
{x_{n + 1}} = {x_0} + \frac{{{h^{q_{n + 1}^{\left( 1 \right)}}}}}{{\Gamma (q_{n
			+ 1}^{(1)} + 2)}}\left( {y_{n + 1}^p + \sum\limits_{j = 0}^n {\rho _{j,n +
			1}^1{y_j}} } \right),\\
{y_{n + 1}} = {y_0} + \frac{{{h^{q_{n + 1}^{(2)}}}}}{{\Gamma (q_{n + 1}^{(2)} +
		2)}}\left( {z_{n + 1}^p + \sum\limits_{j = 0}^n {\rho _{j,n + 1}^2{z_j}} }
\right),\\
{z_{n + 1}} = {z_0} + \frac{{{h^{q_{n + 1}^{(3)}}}}}{{\Gamma (q_{n + 1}^{(3)} +
		2)}}\left( {{f_{n + 1}} - \lambda z_{n + 1}^p - x_{n + 1}^p - {{(x_{n + 1}^p)}^3}
	+ \sum\limits_{j = 0}^n {\rho _{j,n + 1}^3({f_j} - \lambda {z_j} - {x_j} -
		bx_j^3)} } \right),\\
n = 0,1,...,N - 1.
\end{array} \right.
\label{eq8}
\end{equation}
where
\begin{equation}
\left\{ \begin{array}{l}
\rho _{0,n + 1}^i = {n^{q_{_{n + 1}}^{(i)} + 1}} - (n - q_{_{n + 1}}^i){(n +
	1)^{q_{n + 1}^{(i)}}}\\
\rho _{j,n + 1}^i = {(n - j + 2)^{q_{_{n + 1}}^{(i)} + 1}} + {(n - j)^{q_{_{n +
				1}}^{(i)} + 1}} - 2{(n - j + 1)^{q_{_{n + 1}}^{(i)} + 1}}\\
\rho _{n + 1,n + 1}^i = 1,i = 1,2,3
\end{array} \right.
\label{eq9}
\end{equation}
\begin{equation}
\left\{ \begin{array}{l}
x_{n + 1}^p = {x_0} + \frac{{{h^{q_{n + 1}^{(1)}}}}}{{\Gamma (q_{n + 1}^{(1)} +
		1)}}\sum\limits_{j = 0}^n {\theta _{j,n + 1}^1{y_j}} \\
y_{n + 1}^p = {y_0} + \frac{{{h^{q_{n + 1}^{(2)}}}}}{{\Gamma (q_{n + 1}^{(2)} +
		1)}}\sum\limits_{j = 0}^n {\theta _{j,n + 1}^2{z_j}} \\
z_{n + 1}^p = {z_0} + \frac{{{h^{q_{n + 1}^{(3)}}}}}{{\Gamma (q_{n + 1}^{(3)} +
		1)}}\sum\limits_{j = 0}^n {\theta _{j,n + 1}^3({f_j} - \lambda {z_j} - {x_j} -
	bx_j^3)} \\
\theta _{j,n + 1}^i = {(n - j + 1)^{q_{n + 1}^{(i)}}} - {(n - j)^{q_{n +
			1}^{(i)}}},i = 1,2,3.
\end{array} \right.
\label{eq10}
\end{equation}
Here $\Gamma \left(  \cdot  \right)$ is Euler's gamma function. Scheme (10) is a predictor, and scheme (8) is a corrector.

\section{SIMULATION RESULTS}	

\textbf{Test Example.} Let us compare schemes (6) and (8). For this, in the
Cauchy problem (1), we take 
\begin{equation}
f(t) = {t^9} + {t^3} + 6t + \frac{\lambda }{{\Gamma \left( {1 - q(t)}
		\right)}}\frac{d}{{dt}}\left( {\frac{{\Gamma \left( {1 - q(t)} \right)\Gamma
			\left( 4 \right){t^{4 - q(t)}}}}{{\Gamma \left( {5 - q(t)} \right)}}}
\right),{\omega _0} = b = 1,
\label{eq11}
\end{equation}
As a result, we get
\begin{equation}
\begin{array}{l}
\ddot x\left( t \right) + \lambda D_{0t}^{q\left( t \right)}x\left( \tau 
\right) + x\left( t \right) + {x^3}\left( t \right) = {t^9} + {t^3} + 6t +
\frac{\lambda }{{\Gamma \left( {1 - q(t)} \right)}}\frac{d}{{dt}}\left(
{\frac{{\Gamma \left( {1 - q(t)} \right)\Gamma \left( 4 \right){t^{4 -
					q(t)}}}}{{\Gamma \left( {5 - q(t)} \right)}}} \right),\;\\
x\left( 0 \right) = {x_0},\dot x\left( 0 \right) = {y_0},
\end{array}
\label{eq12}
\end{equation}

The exact solution to this problem is the function
\begin{equation}
x(t) = {t^3}.
\label{eq13}
\end{equation}
We will look for the error and computational accuracy by the formulas \cite{7}:
\begin{equation}
{p_i} = \frac{{\ln ({\varepsilon _i})}}{{\ln ({\varepsilon _{i + 1}})}},i =
1,...,N - 1.
\label{eq14}
\end{equation}
In (14), ${x_j}$  is the exact solution, $x_j^{PC}$ is the numerical $\xi_i  = \mathop {\max }\limits_i \left( {\left| {{x_i} - x_i^{PC}} \right|} \right),$ solution, the ${\xi_i}$ -- error at the \textit{i}-th step, the ${\xi _{i + 1}}$ -- error at the \textit{i }+1-th step. Take the following control parameters for the equation (12): $\lambda  = 0.1,\delta  = 0,{x_0} = 0.01,{y_0} = 0.03,T = 1.$
\begin{figure}[h!]
\centering
\includegraphics[scale=1.2]{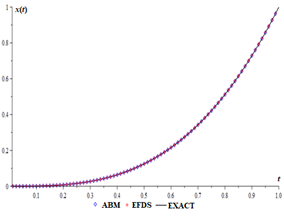}
\caption{The solution of equation (12), obtained according to schemes (6) and (8), as well as the exact solution for $N = 80$}
\end{figure}

\begin{table}[h!]
\caption{Error and computational accuracy for schemes (6) and (8)}
\begin{tabular}{p{21pt}p{49pt}p{49pt}p{49pt}p{56pt}p{49pt}}
	\multicolumn{6}{l}{\parbox{277pt}{
	}} \\
	\parbox{21pt}{\centering } & \parbox{49pt}{\centering } & \multicolumn{2}{l}{\parbox{99pt}{\raggedright 
			\textbf{EFDS (6{\large )}}
	}} & \multicolumn{2}{l}{\parbox{106pt}{\raggedright 
			\textbf{ABM (8{\large )}}
	}} \\
	\hline
	\parbox{21pt}{\centering $N$\textbf{{\small  }}} & \parbox{49pt}{\raggedright \textbf{{\small       $h$ }}} & \parbox{49pt}{\centering $\xi $} & \parbox{49pt}{\centering $p$\textbf{{\small  }}} & \parbox{56pt}{\centering $\xi $} & \parbox{49pt}{\centering $p$} \\
	\hline
	\parbox{21pt}{\centering 
		1{\large 0}
	} & \parbox{49pt}{\centering 
		0.{\large 1}
	} & \parbox{49pt}{\centering 
		{\small 0.011403981}
	} & \parbox{49pt}{\centering 
		-
	} & \parbox{56pt}{\centering 
		0.006415779
	} & \parbox{49pt}{\centering 
		-
	} \\
	\parbox{21pt}{\centering 
		2{\large 0}
	} & \parbox{49pt}{\centering 
		0.0{\large 5}
	} & \parbox{49pt}{\centering 
		{\small 0.008901704}
	} & \parbox{49pt}{\centering 
		{\small 0.947533812}
	} & \parbox{56pt}{\centering 
		0.001892361
	} & \parbox{49pt}{\centering 
		0.805271339
	} \\
	\parbox{21pt}{\centering 
		4{\large 0}
	} & \parbox{49pt}{\centering 
		0.02{\large 5}
	} & \parbox{49pt}{\centering 
		{\small 0.005272949}
	} & \parbox{49pt}{\centering 
		{\small 0.900164648}
	} & \parbox{56pt}{\centering 
		0.000571043
	} & \parbox{49pt}{\centering 
		0.839567683
	} \\
	\parbox{21pt}{\centering 
		8{\large 0}
	} & \parbox{49pt}{\centering 
		0.0125
	} & \parbox{49pt}{\centering 
		{\small 0.002844266}
	} & \parbox{49pt}{\centering 
		{\small 0.894705312}
	} & \parbox{56pt}{\centering 
		0.000173628
	} & \parbox{49pt}{\centering 
		0.86250084
	} \\
	\parbox{21pt}{\centering 
		16{\large 0}
	} & \parbox{49pt}{\centering 
		0.00625
	} & \parbox{49pt}{\centering 
		{\small 0.001474332}
	} & \parbox{49pt}{\centering 
		{\small 0.899210843}
	} & \parbox{56pt}{\centering 
		0.0000525854
	} & \parbox{49pt}{\centering 
		0.878771153
	} \\
	\parbox{21pt}{\centering 
		32{\large 0}
	} & \parbox{49pt}{\centering 
		0.003125
	} & \parbox{49pt}{\centering 
		{\small 0.000750245}
	} & \parbox{49pt}{\centering 
		{\small 0.906108418}
	} & \parbox{56pt}{\centering 
		0.0000158165
	} & \parbox{49pt}{\centering 
		0.891321226
	} \\
	\parbox{21pt}{\centering 
		64{\large 0}
	} & \parbox{49pt}{\centering 
		0.0015625
	} & \parbox{49pt}{\centering 
		{\small 0.000378395}
	} & \parbox{49pt}{\centering 
		{\small 0.913134815}
	} & \parbox{56pt}{\centering 
		0.0000047259
	} & \parbox{49pt}{\centering 
		0.901488248
	} \\
	\parbox{21pt}{\centering 
		128{\large 0}
	} & \parbox{49pt}{\centering 
		0.00078125
	} & \parbox{49pt}{\centering 
		{\small 0.000190019}
	} & \parbox{49pt}{\centering 
		{\small 0.919609734}
	} & \parbox{56pt}{\centering 
		0.0000015483
	} & \parbox{49pt}{\centering 
		0.9165891
	} \\
	\hline
\end{tabular}
\end{table}

\textbf{Example.} Let's take the following parameters for the system (3):

\[
\lambda=\delta=\omega=\omega_0=b=1, x_0=y_0=0, T=100, N=1800, q\left(t\right)=0.8-\frac{t}{2T}.
\]

\begin{figure}[h!]
	\begin{minipage}[h!]{0.44\linewidth}
		\center{\includegraphics[width=0.8\linewidth]{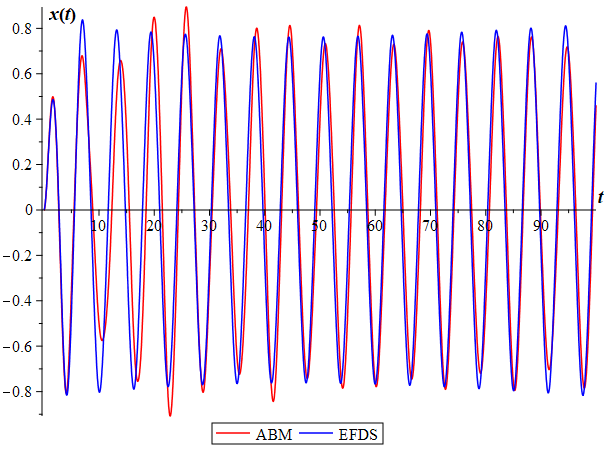} \\ (a)}
	\end{minipage}
	\hfill
	\begin{minipage}[h!]{0.44\linewidth}
		\center{\includegraphics[width=0.8\linewidth]{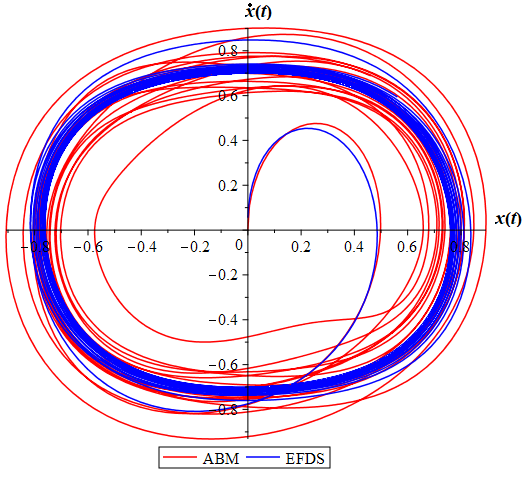} \\ (b)}
	\end{minipage}
	\caption{Phase trajectory (a) and oscillogram (b) for
		numerical schemes (6) and (8)}
	\label{ris:image1}
\end{figure}

Figure 2 shows the phase trajectory (a) and the oscillogram (b) for the
finite-difference scheme (6) and the predictor-corrector (8). Oscillograms
characterize steady fluctuations and tend to a regular mode. The phase
trajectories characterize the limit cycle, which is determined by an external
harmonic action with an amplitude $\delta $ and a frequency $\omega $.

\begin{table}[h!]
\centering
\caption{Error and computational accuracy for schemes (6) and (8)}
\begin{tabular}{p{21pt}p{49pt}p{49pt}p{49pt}p{56pt}p{49pt}}
	\multicolumn{6}{l}{\parbox{277pt}{
	}} \\
	\parbox{21pt}{\centering } & \parbox{49pt}{\centering } & \multicolumn{2}{l}{\parbox{99pt}{\raggedright 
			\textbf{ABM (8{\large )}}
	}} & \multicolumn{2}{l}{\parbox{106pt}{\raggedright 
			\textbf{EFDS (6{\large )}}
	}} \\
	\hline
	\parbox{21pt}{\centering $N$\textbf{{\small  }}} & \parbox{49pt}{\raggedright \textbf{{\small       $h$ }}} & \parbox{49pt}{\centering $\xi $} & \parbox{49pt}{\centering $p$\textbf{{\small  }}} & \parbox{56pt}{\centering $\xi $} & \parbox{49pt}{\centering $p$} \\
	\hline
	\parbox{21pt}{\centering 
		1{\large 0}
	} & \parbox{49pt}{\centering 
		0.{\large 1}
	} & \parbox{49pt}{\centering 
		{\small 0.023757}
	} & \parbox{49pt}{\centering 
		-
	} & \parbox{56pt}{\raggedleft 
		{\small 0.012298}
	} & \parbox{49pt}{\centering 
		-
	} \\
	\parbox{21pt}{\centering 
		2{\large 0}
	} & \parbox{49pt}{\centering 
		0.0{\large 5}
	} & \parbox{49pt}{\centering 
		{\small 0.00845}
	} & \parbox{49pt}{\raggedleft 
		{\small 0.783443}
	} & \parbox{56pt}{\raggedleft 
		{\small 0.004977}
	} & \parbox{49pt}{\raggedleft 
		{\small 0.829423}
	} \\
	\parbox{21pt}{\centering 
		4{\large 0}
	} & \parbox{49pt}{\centering 
		0.02{\large 5}
	} & \parbox{49pt}{\centering 
		{\small 0.003635}
	} & \parbox{49pt}{\raggedleft 
		{\small 0.849825}
	} & \parbox{56pt}{\raggedleft 
		{\small 0.004018}
	} & \parbox{49pt}{\raggedleft 
		{\small 0.961204}
	} \\
	\parbox{21pt}{\centering 
		8{\large 0}
	} & \parbox{49pt}{\centering 
		0.0125
	} & \parbox{49pt}{\centering 
		{\small 0.001768}
	} & \parbox{49pt}{\raggedleft 
		{\small 0.886322}
	} & \parbox{56pt}{\raggedleft 
		{\small 0.011872}
	} & \parbox{49pt}{\raggedleft 
		{\small 1.244357}
	} \\
	\parbox{21pt}{\centering 
		16{\large 0}
	} & \parbox{49pt}{\centering 
		0.00625
	} & \parbox{49pt}{\centering 
		{\small 0.000905}
	} & \parbox{49pt}{\raggedleft 
		{\small 0.904421}
	} & \parbox{56pt}{\raggedleft 
		{\small 0.020541}
	} & \parbox{49pt}{\raggedleft 
		{\small 1.141099}
	} \\
	\parbox{21pt}{\centering 
		32{\large 0}
	} & \parbox{49pt}{\centering 
		0.003125
	} & \parbox{49pt}{\centering 
		{\small 0.000469}
	} & \parbox{49pt}{\raggedleft 
		{\small 0.914243}
	} & \parbox{56pt}{\raggedleft 
		{\small 0.030448}
	} & \parbox{49pt}{\raggedleft 
		{\small 1.112723}
	} \\
	\hline
\end{tabular}
\end{table}

Table 2 presents an estimate of the computational accuracy according to schemes (6) and (8), obtained by the double recalculation method (Runge's rule). We see that with an increase in the nodes of the computational grid, the computational accuracy tends to unity. It should also be noted that for this example the predictor-corrector method converges faster than the explicit finite-difference scheme.

Note that the computational accuracy for an explicit finite-difference scheme
does not behave stably, but nevertheless tends to unity. It is necessary to carry out more rigorous estimates of the stability and convergence of the method.

\section{Conclusion}

In this paper, we have proposed two numerical methods for solving the
mathematical model of the Duffing fractional oscillator. The first method is
based on the approximation of the derivative of a fractional variable order based on the Grunwald-Letnikov difference operator, and the second is the
Adams-Bashfort-Moulton method from the class of predictor-corrector methods. It is shown by examples that the ABM method has better convergence. Further
continuation of the work is the study of the stability and convergence of the
method of the explicit finite-difference scheme.

\begin{acknowledgments}
The work was performed within the framework of the research project of Vitus
Bering Kamchatka State University "Natural disasters in Kamchatka - earthquakes and volcanic eruptions (monitoring, forecast, study, psychological support of the population)" no. AAAA-A19-119072290002-9.
\end{acknowledgments}


\end{document}